\documentclass [12pt] {article}
\usepackage {amsmath}
\usepackage[latin1]{inputenc}
\usepackage {amsfonts, amssymb,amsthm}
\usepackage {latexsym}
\newtheorem {theorem} {Theorem} [section]
\newtheorem {lemma} {Lemma} [section]
\newtheorem {cor} {Corollary} [section]
 \newtheorem{preremark}{Remark}[section]
    {\begin{preremark}\rm}{\end{preremark}}
     \newtheorem{preremark1}{Example}[section]
    {\begin{preremark1}\rm}{\end{preremark1}}
         \newtheorem{preremark2}{Definition}[section]
  \newenvironment{defn}%
    {\begin{preremark2}\rm}{\end{preremark2}}

\begin{document}
\title{The second boundary value problem for equations of viscoelastic diffusion in polymers}

\author{Dmitry A. Vorotnikov}

\date{}

\maketitle

\begin{center}

CMUC (Centro de Matemática da Universidade de Coimbra)

Apartado 3008, 3001 - 454 Coimbra, Portugal

mitvorot@mat.uc.pt

\end{center}

\begin{abstract} The classical approach to diffusion processes is based on Fick's law that the flux is proportional to the concentration gradient. Various phenomena occurring during propagation of penetrating liquids in polymers show that this type of diffusion exhibits anomalous behavior and contradicts the just mentioned law. However, they can be explained in the framework of non-Fickian diffusion theories based on viscoelasticity of polymers.  Initial-boundary value problems for viscoelastic diffusion equations have been studied by several authors. Most of the studies are devoted to the Dirichlet BVP (the concentration is given on the boundary of the domain). In this chapter we study the second BVP, i.e. when the normal component of the concentration flux is prescribed on the boundary, which is more realistic in many physical situations. We establish existence of weak solutions to this problem.  We suggest some conditions on the coefficients and boundary data under which all the solutions tend to the homogeneous state as time goes to infinity.

\end{abstract}

\section {Introduction}

\renewcommand{\theequation}{\arabic{section}.\arabic{equation}}

\numberwithin{equation}{section}

\newcommand {\R} {\mathbb{R}}
\newcommand {\E} {\mathbf{E}}
\def\be{\begin{equation}}
\def\ee{\end{equation}}
\def\fr#1#2{\frac{\partial #1}{\partial #2}}

The continuity equation for diffusion \be\fr u t = - div\, J\ee states that variations of the concentration $u(t,x)$ at any spatial
point $x$ and moment of time $t$ can only be caused by inflow and outflow of a penetrant into and out of that area. Here $J=J(t,x)$ is the
concentration flux vector.

The classical diffusion theory is based on Fick's law (the flux is proportional to
the concentration gradient with negative proportionality factor $-D$). The continuity equation and Fick's law yield the
classical diffusion equation \be\fr u t = div\, (D(u) \nabla
u),\ee
which becomes the heat equation \be\fr u t
= D \Delta u\ee for constant diffusion coefficient $D$.

The concentration behaviour for diffusion of penetrant liquids in
poly\-mers exhibits such phenomena as \textit{case II diffusion},
\textit{sorption overshoot}, \textit{literal skinning, trapping
skinning} and \textit{desorption overshoot}, which contradict
Fick's law, see e.g. \cite{chn1, ed,ed2, cc, tw,tw1, var, wit}.
There is a number of approaches which explain these non-Fickian
properties of polymeric diffusion. They are usually based on
taking into account the viscoelastic nature of polymers (cf.
\cite{lee} and references therein) and on the possibility of
glass-rubber phase transition (see e.g. \cite{wit} with some
review). We are going to study the model which is due to Cohen et
al. (see \cite{chn0,chn1,chn3}; related models or particular cases
of this one were suggested by other authors, see e.g.
\cite{dur,wit}). It consists in combining the continuity equation
with the following system

\be J(t,x)= -D\nabla u - E \nabla \sigma+M u, \ee
\be \fr
{\sigma} t + \beta  \sigma = \mu u+ \nu u'. \ee

Generally speaking, the coefficients $\beta$, $D$, $E$, $M$, $\mu$ and $\nu$ may depend $u$, or even on $t,x$ and $\sigma$.

Let us briefly discuss the meaning and typical behaviour of the coefficients. The scalar
function $\beta$ is the inverse of the relaxation time. A characteristic
form of $\beta$ is \cite{chn1} \be\beta(u)= \frac 1 2 (\beta_R+
\beta_G)+ \frac 1 2 (\beta_R- \beta_G) \tanh (\frac {u- u_{RG}}
{\delta})\ee where $\beta_R,\beta_G, \delta, u_{RG}$ are positive
constants, $\beta_R>\beta_G$. The polymer-penetrant systems modeled with the help of (1.6) can be in two phases: glassy and rubbery. The glassy state corresponds to the areas of low concentration. Here the polymer network is severely entangled, and the relaxation time is high, so its inverse is low. Moreover, it is close to a certain value $\beta_G$. In the high concentration areas the system is in the rubbery state: the network disentangles, so the relaxation time is small, and its inverse is close to $\beta_R>\beta_G$. The glass-rubber phase transition occurs near a certain concentration $u_{RG}$, and the value of $\delta$ determines the length of the transition segment.
The coefficients $D$ and $E$ are non-negative scalars (more generally, they are positive-definite tensors) called
the diffusion and stress-diffusion coefficients, respectively. As the concentration increases, the polymer network disentangles, so the diffusivity also increases. Thus, $D$ should be an increasing function of concentration: in particular, $D$ can depend on $u$ in a way similar to (1.6) \cite{cc}. $E$ is sometimes considered to be a constant, see \cite{chn0} for some justification, but numerical simulations \cite{cox} have shown that, if $E(0)\neq 0$, then the concentration $u$ may become negative, which is physically meaningless. Conversely, it can analytically be proved  that, if $E(0)= 0$, then the concentration $u$ remains non-negative provided it is non-negative at the initial moment of time \cite{am2}. In \cite{iter}, we make related observations showing expediency of the condition $E(1)= 0$, which can maintain the concentration $u$ of less than or equal to $100\%$. Thus, a modeling example is \be E(u)=\frac {\alpha_1 u(u-1)^2}{\alpha_2 +(u-1)^2},\ee where $\alpha_1$ and $\alpha_2$ are positive constants, $\alpha_2$ is small.
The functions $\mu$ and $\nu$ should be non-negative and bounded \cite{chn0}, $M$ is the convection velocity vector, assumptions on it will be given below.

 Let $\mathrm{n}(x)$ be the outward normal vector at the point $x$ of the boundary $\partial\Omega$ of a domain $\Omega\subset\mathbb{R}^n$ \footnote{The most important particular cases are $n=2$
(diffusion in polymer films) and $n=3$.}. Then system (1.1),(1.4),(1.5) may be completed with such boundary conditions as
\be u(t,x)=\phi(t,x),\ x\in\partial\Omega\ee (the concentration on the boundary is prescribed) and \be -\sum\limits_{i=1}^n J_i(t,x) \mathrm{n}_i(x)=\varphi(t,x),\ x\in\partial\Omega\ee (the influx\footnote{It can be negative.} of the penetrant through the boundary is known).

The initial-boundary value problems for system (1.1),(1.4),(1.5) possess maximal (not global
in time) solutions for a more general boundary condition, which includes (1.8) and (1.9), see \cite{am2}.
The global (in time) existence results are known for the Dirichlet condition (1.8). A theorem on global  solvability is presented in
\cite{am1} for $f=\mu u$, $M\equiv 0$ and $D=E$ being a constant scalar. It is formulated for
the one-dimensional case ($0<x<1$), but the technique used there seems to be applicable for
$x\in \Omega $, where $\Omega\subset \mathbb{R}^n$ is a bounded
domain with a smooth boundary. Another global existence result is given in
\cite{bei}. They assume the stress-diffusion coefficient $E$ to be a non-constant increasing function of concentration,  $E(0)=0$ (so $E(1)=0$ is not allowed). However, it is required that the initial and boundary data for the concentration are bounded from below by a positive constant, so the solution is always strictly positive, and this approach does not permit to consider dry regions in a polymer. Paper \cite{bei} is mainly concerned with the one-dimensional case, but also they suggest a brief plan how to generalize the result for the multidimensional situation. Global existence of dissipative (ultra weak) solutions for constant scalar $D$ and $E$ and $M\equiv 0$ is shown in \cite{diss} for $\Omega=\mathbb{R}^n$ (again, the ideas used there seem to be suitable for
$\Omega\subset \mathbb{R}^n$). In \cite{var}, global (in time) weak solutions on a bounded domain $\Omega \subset \R^n$ are constructed, under rather general assumptions on the coefficients. Further investigation of the weak solutions of the Dirichlet problem is carried out in \cite{iter}: it is proved that, for any sufficiently short time segment and any stress prescribed at the beginning of this segment, there exists a weak solution such that the concentrations at the beginning and at the end of the segment are the same, and, under an additional assumption on coefficients, existence of time-periodic weak solutions (without any restrictions of the period length) is shown. Paper \cite{sell} considers long-time behaviour issues for this problem: provided $D$ and $E$ are constant scalars and $M\equiv 0$, the solutions generate a dissipative semiflow, and there exist a minimal trajectory attractor and a global attractor.

In this chapter, we construct (global in time) weak solutions for problem (1.1),(1.4),(1.5),(1.9) on a bounded domain $\Omega \subset \R^n$ for given initial concentration and stress.
The coefficients may depend on $t,x,u$ and $\sigma$. In addition, we suggest some conditions on the coefficients and boundary data under which all the solutions tend to the homogeneous state $u=const$ as time goes to infinity. The chapter
is organized
in the following way. In Section 2, we introduce the required notations.
In Section 3, we give a weak formulation of the initial-boundary value problem and state the result on existence of weak solutions (Theorem 3.1), which is proved in Section 4. In Section 5, we touch the long-time behaviour.

\section {Notation}
We use the standard notations $L_p (\Omega) $, $W_p^{m} (\Omega)
$, $H^{m} (\Omega) =$ $W_2^{m} (\Omega) $ for Lebesgue and Sobolev
spaces of functions defined on a bounded open set (domain) $\Omega\subset
\R^n$, $n\in \mathbb{N}$.

The scalar product and the Euclidean norm in
$L_2(\Omega)^k=L_2(\Omega, \R^k)$ are denoted by $(u,v)$ and
$\|u\|$, respectively ($k$ is equal to $1$ or $n$). In $H^m(\Omega)$, $m\in \mathbb{N}$, we use the scalar product
$(u,v)_m=\sum\limits_{|\alpha|\leq m}(D^\alpha u, D^\alpha v)$ and the corresponding Euclidean norm $\|u\|_m$.

The space of linear continuous functionals on
$H^m(\Omega)$ (the dual space) is denoted by $H^{-m}_N(\Omega)$.
The value of a functional from
$H ^ {-m}_N(\Omega) $ on an element from $H^m(\Omega) $ is denoted
by $ \langle\cdot, \cdot\rangle $ (the "bra-ket" notation). Similarly, the dual space of $W^m_p(\Omega)$  is denoted by $W^{-m}_{q,N}(\Omega)$, $\frac 1 p +\frac 1 q=1$, $1<p<\infty$, with the corresponding use of the "bra-ket" notation.

Sometimes we shall write simply $L_p $, $H^m$ for $L_p (\Omega)^k, H^m (\Omega)^k$
etc., $k=1,n$.

Let us introduce some basic operators. The operator $div_N:L_q(\Omega)\to W^{-1}_{q,N}(\Omega)$ is determined by the formula \be\langle div_N v, \phi\rangle=-\int\limits_{\Omega}v(x)\nabla \phi(x)\; dx,\ \phi\in W^1_p(\Omega), \frac 1 p +\frac 1 q=1.\ee

The isomorphic operators $A: H^1(\Omega)\to H^{-1}_N(\Omega)$ and $A_2: H^2(\Omega)\to H^{-2}_N(\Omega)$ are given by the expressions \be\langle A v, \phi\rangle=(v,\phi)_1,\ \phi\in H^1(\Omega), \ee \be\langle A_2 v, \phi\rangle=(v,\phi)_2,\ \phi\in H^2(\Omega).\ee
Note that $L_2(\Omega)\subset H^{-1}_N(\Omega)\subset H^{-2}_N(\Omega)$ with natural imbedding operators, and then $Av=v-div_N \nabla v$, $v\in H^1$.

Set $X_N=X_N(\Omega)=A^{-1}(H^{1}(\Omega))$. The scalar product and norm in $X_N$ are $(u,v)_X = (A u, A v)_1$, $\|u\|_X = \|A u\|_1$. The duality between $H^{-1}_N(\Omega)$ and $X_N(\Omega)$ is given by the formula
\be\left\langle u,v\right\rangle _1=\left\langle u,A v\right\rangle,\ u\in H^{-1}_N, v\in X_N.\ee
Note that $\left\langle u,v\right\rangle _1=(u,v)_1$ for $ u\in H^1, v\in X_N$.
The elements of $X_N$ are solutions of the Neumann problem
\be v-\Delta v=u\in H^1(\Omega),\ee
\be\fr{v} {\mathrm{n}} (x)=0, x\in \partial \Omega.\ee
Thus, \be X_N(\Omega)\subset H^2(\Omega)\subset W^1_{2n/n-2}(\Omega)\ee by Sobolev theorem (for sufficiently regular $\Omega$).

The symbols $C (\mathcal{J}; E) $, $C_w (\mathcal{J}; E) $, $L_2
(\mathcal{J}; E) $ etc. denote the spaces of continu\-ous, weakly
continuous, quadratically integrable etc. functions on an interval
$\mathcal{J}\subset \mathbb {R} $ with
values in a Banach space $E $. We recall that a function $u:
\mathcal{J} \rightarrow E$ is \textit{weakly continuous} if for
any linear continuous functional $g$ on $E$ the function  $g(
u(\cdot)): \mathcal{J}\to \mathbb{R}$ is continuous.

If $E$ is a function space ($L_2(\Omega), H^m(\Omega)$ etc.), then
we identify the elements of $C (\mathcal{J}; E) $,
$L_2(\mathcal{J};E)$ etc. with scalar functions defined on
$\mathcal{J}\times \Omega$ according to the formula
$$u(t)(x)=u(t,x),\, t\in \mathcal{J}, x\in \Omega.$$

We shall also use the function spaces ($T$ is a positive number):
$$ W_N= W_N(\Omega, T) = \{\tau\in L_2 (0, T; H^1 (\Omega)), \
\tau' \in L_2 (0, T;
 H ^ {-1}_N (\Omega)) \}$$ $$\|\tau\|_{W_N}=\|\tau\|_{L_2 (0, T; H^1(\Omega))}+\|\tau'\|_{L_2 (0, T; H ^ {-1}_N(\Omega))};$$

$$ W_1= W_1(\Omega, T) = \{\tau\in L_2 (0, T; X _N(\Omega)), \
\tau' \in L_2 (0, T;
 H ^ {-1}_N (\Omega)) \}$$ $$\|\tau\|_{W_1}=\|\tau\|_{L_2 (0, T; X_N(\Omega))}+\|\tau'\|_{L_2 (0, T; H ^ {-1}_N(\Omega))};$$

 $$ W_2= W_2(\Omega, T) = \{\tau\in L_2 (0, T; H^2 (\Omega)), \
\tau' \in L_2 (0, T;
 H ^ {-2}_N (\Omega)) \}$$ $$\|\tau\|_{W_2}=\|\tau\|_{L_2 (0, T; H^2(\Omega))}+\|\tau'\|_{L_2 (0, T; H ^ {-2}_N(\Omega))}.$$

Lemma III.1.2 from \cite{tem} implies continuous embeddings $W_N,W_2\subset
C([0,T];L_2(\Omega))$, $W_1\subset
C([0,T];H^1(\Omega))$ (see also \cite{ggz}).

We use the notation $|\cdot|$ for the absolute value of a number, for the Euclidean
norm in $\R^n$, and in the following case.

Denote
by $\R^{n \times n}$ the space of  matrices of the order $n\times n$ with the norm
$$|Q|=\max_{\xi\in \R^n,|\xi|=1} |Q\xi|.$$
Let $\R^{n \times n}_+ \subset$ $\R^{n \times n}$ be the set of such matrices $Q$ that
$$ (Q\xi, \xi)_{\R^n} \geq d(Q)(\xi, \xi)_{\R^n}$$ for some $d(Q)\geq 0$ and all $\xi\in \R^n$.

The symbol $C$ will stand for a generic positive constant that
can take different values in different places.

\section {Weak formulation of the problem}
We consider a polymer filling a sufficiently regular\footnote{Say,
it is locally located on one side of its $C^2$-smooth boundary.}
bounded domain $\Omega\subset\mathbb{R}^n$, $n\in \mathbb{N}$. We study the diffusion of a penetrant in this
polymer which is described\footnote{This problem makes sense for
diffusion in polymers provided $0\leq u\leq 1$, i.e. the
concentration is not less than $0\%$ and does not exceed $100\%$.
The assumptions on coefficients which guarantee this condition
provided it is fulfilled at the initial moment of time are
discussed in \cite{iter} for the Dirichlet problem, mainly based
on the results of \cite{am2}; similar arguments are applicable in
the Neumann case (1.9) with $\varphi\equiv 0$. However, this question is
still not completely studied. So we consider here the general
setting (3.1) -- (3.4).} by the following initial-boundary value
problem\footnote{System (3.1), (3.2), (3.3) is obtained from
(1.1),(1.4),(1.5),(1.9). For technical purposes, we assign
subscript zero to the coefficients.}:

\be \fr u t = div[D_0(t,x,u,\sigma) \nabla u $$ $$+ E_0(t,x,u,\sigma) \nabla \sigma -M_0(t,x,u,\sigma) u],\ (t,x)\in [0,T]\times \Omega,\ee \be \fr {\sigma} t
+ \beta_0 (t,x,u,\sigma) \sigma = \mu_0(u) u + \nu_0(u) \fr u t,\ (t,x)\in [0,T]\times \Omega, \ee \be \sum\limits_{i,j=1}^n \Big [ D_0(t,x,u,\sigma)_{ij} \fr u {x_j}+ E_0(t,x,u,\sigma)_{ij} \fr \sigma {x_j} $$ $$-M_0(t,x,u,\sigma)_i u \Big ] \mathbf{n}_i(x)=\varphi(t,x),\ (t,x)\in [0,T]\times \partial\Omega, \ee \be u(0,x)
=u_0(x), \ \sigma(0,x) = \sigma_0(x),\ x\in\Omega.\ee

Here $u=u(t,x):[0,T]\times \overline{\Omega} \to \R$ is the unknown
concentration of the penetrant (at the spatial point $x$ at the
moment of time $t$), $\sigma=\sigma(t,x):[0,T]\times \overline{\Omega}
\to \R$ is the unknown stress, $u_0=u_0(x)$, $\sigma_0=\sigma_0(x):
\Omega \to \R$ are given initial data, $\varphi:[0,T]\times \partial\Omega \to \R$ is the influx of the liquid through the boundary, $\mu_0,\nu_0:\R\to \R$, $D_0$, $E_0: \R^{n+3}=\R\times\R^n\times\R\times\R\to \R^{n\times n}_+,\ \ $ $\beta_0: \R^{n+3}\to \R$, $M_0: \R^{n+3}\to \R^n$ are given functions, $\nu_0(\cdot)\geq 0$.

Before giving a definition of a weak solution to this problem, it is convenient to make a change of variables. Denote $$\varsigma(t,x)=\sigma(t,x)-\int\limits_0^{u(t,x)}\nu_0(y)\;dy,$$ $$ \varsigma_0(x)=\sigma_0(x)-\int\limits_0^{u_0(x)}\nu_0(y)\;dy,$$
 $$D(t,x,u,\varsigma)=$$ $$D_0\left(t,x,u,\varsigma+\int\limits_0^u \nu_0(y)dy\right)+\nu_0(u)E_0\left(t,x,u,\varsigma+\int\limits_0^u \nu_0(y)dy\right)\in \R^{n\times n}_+ ,$$ $$E(t,x,u,\varsigma)=E_0\left(t,x,u,\varsigma+\int\limits_0^u \nu_0(y)dy\right),$$ $$f(t,x,u,\varsigma)=-u M_0\left(t,x,u,\varsigma+\int\limits_0^u \nu_0(y)dy\right),$$ $$\beta_1(t,x,u,\varsigma)=-\beta_0\left(t,x,u,\varsigma+\int\limits_0^u \nu_0(y)dy\right),$$ $$\gamma\left(t,x,u,\varsigma\right)=\mu_0(u)-\frac{\beta_0\left(t,x,u,\varsigma+\int\limits_0^u \nu_0(y)dy\right)\int\limits_0^{u}\nu_0(y)\;dy }{u}.$$
Note that, if $u$ vanishes, then, by continuity, we consider the last term to become

$$-\beta_0(t,x,0,\varsigma)\nu_0(0).$$
Then we
can rewrite (3.1) -- (3.4) in the following form:

\be \fr u t = div[D(t,x,u,\varsigma) \nabla u + E(t,x,u,\varsigma) \nabla \varsigma + f(t,x,u,\varsigma) ], \ee \be \fr {\varsigma} t =
\beta_1 (t,x,u,\varsigma) \varsigma + \gamma(t,x,u,\varsigma) u, \ee
\be \sum\limits_{i,j=1}^n \Big [ D(t,x,u,\varsigma)_{ij} \fr u {x_j}+ E(t,x,u,\varsigma)_{ij} \fr \varsigma {x_j} $$ $$+f(t,x,u,\varsigma)_i  \Big ] \mathbf{n}_i(x)=\varphi(t,x),\ (t,x)\in [0,T]\times \partial\Omega, \ee
 \be u | _ {t=0} =u_0, \ \varsigma| _
{t=0} = \varsigma_0.\ee

Now, before describing our assumptions on the coefficients, let us calculate the gradient of the right member of (3.6):
\be\nabla
(\beta_1 (t,x,u,\varsigma) \varsigma) + \nabla(\gamma(t,x,u,\varsigma) u)$$ $$ =
\beta_1 (t,x,u,\varsigma) \nabla\varsigma + \fr{\beta_1}{x} (t,x,u,\varsigma) \varsigma +\fr{\beta_1}{u} (t,x,u,\varsigma) \varsigma\nabla u +\fr{\beta_1}{\varsigma} (t,x,u,\varsigma) \varsigma\nabla\varsigma $$ $$ + \gamma(t,x,u,\varsigma) \nabla u+\fr \gamma x (t,x,u,\varsigma) u +\fr \gamma u (t,x,u,\varsigma) u \nabla u +\fr \gamma \varsigma (t,x,u,\varsigma) u \nabla\varsigma $$ $$=\beta(t,x,u,\varsigma) \nabla u+ \mu(t,x,u,\varsigma) \nabla \varsigma +g(t,x,u,\varsigma), \ee where
\be\beta(t,x,u,\varsigma)=\fr{\beta_1}{u} (t,x,u,\varsigma)\varsigma
+ \gamma(t,x,u,\varsigma)+\fr \gamma u
(t,x,u,\varsigma) u, \ee
\be\mu(t,x,u,\varsigma)=\beta_1 (t,x,u,\varsigma)  +
\fr{\beta_1}{\varsigma} (t,x,u,\varsigma) \varsigma
+\fr \gamma \varsigma (t,x,u,\varsigma) u,\ee
\be g(t,x,u,\varsigma)=  \fr{\beta_1}{x} (t,x,u,\varsigma) \varsigma+\fr \gamma x (t,x,u,\varsigma) u.\ee


We assume the following:
\\
i) $D$, $E: \R^{n+3}\to \R^{n\times n};$ $f,g: \R^{n+3}\to \R^n$; $\mu$, $\beta,\gamma,\beta_1: \R^{n+3}\to \R$.\\
ii) Each of these eight functions (e.g. $D(t,x,u,\varsigma)$) is measurable in
 $(t,x)$ for fixed $(u,\varsigma)$.\\
iii) Each of these functions is continuous in
 $(u,\varsigma)$ for fixed $(t,x)$.\\
iv) These functions satisfy the estimates
\be |D(t,x,u,\varsigma)|\leq K_D,\ee
\be |E(t,x,u,\varsigma)|\leq K_E,\ee
\be \max(|\beta(t,x,u,\varsigma)|,  |\gamma(t,x,u,\varsigma)|)\leq K_\beta,\ee
\be \max(|\mu(t,x,u,\varsigma)|, |\beta_1(t,x,u,\varsigma)|)\leq K_\mu,\ee
\be |f(t,x,u,\varsigma)|\leq K_f(|u|+|\varsigma|)+ \widetilde f(t,x),\ee
\be |g(t,x,u,\varsigma)|\leq K_g(|u|+|\varsigma|)+ \widetilde g(t,x)\ee
with some constants $K_D,\dots,K_g$ and functions\footnote{Clearly, the behaviour of these functions outside $(0,T)\times\Omega$ does not matter.} $\widetilde f, \widetilde g \in L_{2,loc}(\R^{n+1})$.\\
v) \be (D(t,x,u,\varsigma)\xi, \xi)_{\R^n} \geq d(\xi, \xi)_{\R^n},\ee
where $d>0$ is independent of $(t,x,u,\varsigma)\in \R^{n+3}$ and $\xi\in  \R^{n}$.\\
vi) Relations (3.10) -- (3.12) hold.

It is easy to see that, if $E_0$ and $\beta_0$ are taken in the forms (1.7) and (1.6), then (3.14) and (3.15) are violated. It turns out that such deficiencies can be corrected without loss of generality of the model (see \cite[Section 3]{iter} for a detailed discussion, cf. also \cite{var,sell}). In brief, physically, the concentration $u$ and the stress $\varsigma$ are uniformly bounded,
therefore the coefficients of systems (3.1)--(3.2) (and, consequently, of (3.5)--(3.6)) can be experimentally determined only for bounded $u$ and $\varsigma$, whereas "at infinity" we can choose them at discretion.

Let us now rewrite (3.5) and (3.7) in a weak form. Assuming $u$ and $\varsigma$ sufficiently regular, take the $L_2(\Omega)$-scalar product of the members of (3.5) with a test function $\phi\in H^1(\Omega)$, and integrate by parts in the right-hand side:

\be (u',\phi)= -(D(t,x,u,\varsigma) \nabla u + E(t,x,u,\varsigma) \nabla \varsigma + f(t,x,u,\varsigma),\nabla \phi)$$ $$ + \sum\limits_{i,j=1}^n \int\limits_{\partial\Omega}\Big [ D(t,x,u,\varsigma)_{ij} \fr u {x_j}+ E(t,x,u,\varsigma)_{ij} \fr \varsigma {x_j} +f(t,x,u,\varsigma)_i  \Big ] \mathbf{n}_i(x)\phi(x)\,ds $$ $$= -(D(t,x,u,\varsigma) \nabla u + E(t,x,u,\varsigma) \nabla \varsigma + f(t,x,u,\varsigma),\nabla \phi) + \int\limits_{\partial\Omega}\varphi\phi\,ds. \ee
Denote by $\psi(t)$ the linear functional $\phi\mapsto \int\limits_{\partial\Omega}\varphi(t)\phi\,ds$. We assume that (for a.a. $t$) this integral exists and continuously depends on $\phi\in H^1(\Omega)$, so $\psi(t)\in H^{-1}_N(\Omega)$; clearly, this is true e.g. if $\varphi(t)\in L_2(\partial \Omega)$. Then we arrive at
\be \fr u t = div_N[D(t,x,u,\varsigma) \nabla u + E(t,x,u,\varsigma) \nabla \varsigma + f(t,x,u,\varsigma) ] +\psi,\ee
which should be understood as an equality of functionals from $H^{-1}_N(\Omega)$. Conversely, for each pair of sufficiently regular functions $(u,\varsigma)$, (3.21) implies (3.5) and (3.7).

\begin{defn} \rm A pair of
functions $(u,\varsigma)$ from the class \be u \in W_N(\Omega,T), \varsigma\in
H^1(0,T;H^1(\Omega))\ee is a {\it weak} solution to problem (3.5)-(3.8)
if it satisfies (3.8), equality (3.21) holds in the space $H^{-1}_N(\Omega)$ a.e. on
$(0,T)$, and (3.6) holds a.e. in
$(0,T)\times \Omega$. \end{defn}

Note that (3.8) makes sense due to the embeddings $$W_N\subset C([0,T];L_2(\Omega)),\ H^1(0,T;H^1(\Omega))\subset
C([0,T];H^1(\Omega)).$$

\begin{theorem} For every $u_0 \in L_2(\Omega)$,  $\varsigma_0\in H^1(\Omega)$ and $\psi\in L_{2}(0,T;H^{-1}_N(\Omega))$, there exists a weak solution to problem (3.5) -- (3.8) in class (3.22). \end{theorem}

\section {Proof of the existence result}

The proof of Theorem 3.1 is based on the study of
the following auxiliary problem:

\be \fr v t + \varepsilon A_2 v= \lambda div_N[D(t,x,v,\tau) \nabla v + E(t,x,v,\tau) \nabla \tau +f(t,x,v,\tau)]+\lambda\psi, \ee \be \fr {\tau} t  + \varepsilon A^2\tau=\lambda[
\beta_1 (t,x,v,\tau) \tau + \gamma(t,x,v,\tau) v], \ee  \be v | _ {t=0} =u_0,\ee \be \tau| _ {t=0} =
\varsigma_0.\ee Here $\varepsilon>0$, $\lambda\in[0,1]$ are parameters.
We are going to derive some a priori estimates for the weak solutions of this problem. Then we shall show its solvability via topological degree arguments (the presence of the parameter $\lambda$ is important at this stage). Finally, we shall put $\lambda=1$ and pass to the limit as $\varepsilon \to 0$.

\begin{defn} \rm Given $u_0 \in L_2(\Omega)$, $\varsigma_0\in H^1(\Omega)$, a pair of
functions $(v,\tau)$ from the class \be v \in W_2(\Omega,T),\tau\in
W_1(\Omega,T)\ee is a {\it weak} solution of problem (4.1)-(4.4)
if equality (4.1) holds in the space $H^{-2}_N(\Omega)$ a.e. on
$(0,T)$, (4.2) holds in the space $H^{-1}_N(\Omega)$ a.e. on
$(0,T)$, (4.3) holds in $L_2(\Omega)$, and (4.4) holds in $H^1(\Omega)$.\end{defn}

The last two conditions make sense due to the embeddings $$W_1\subset C([0,T];H^1(\Omega)),\ W_2\subset
C([0,T];L_2(\Omega)).$$

\begin{lemma} Let $(v,\tau)$ be a weak solution to problem
(4.1)-(4.4). Then the following a priori estimate holds: \be
\varepsilon\|v\|^2_{L_2(0,T;H^2(\Omega))}+ \varepsilon\|\tau\|^2_{L_2(0,T;X_N)}+$$ $$ \|v\|^2_{L_\infty(0,T;L_2(\Omega))} + \lambda\|v\|^2_{L_2(0,T;H^1(\Omega))} +
\|\tau\|^2_{L_\infty(0,T;H^1(\Omega))}\leq C \ee where $C$ is
independent of $\lambda$ and $\varepsilon$. \end{lemma}
\textbf{Proof.} Take the "bra-ket" of the terms of (4.2) (as
elements of $H^{-1}_N(\Omega)$) and $A\tau(t)\in H^1(\Omega)$ at a.a.
$t\in[0,T]$: \be \left\langle  \tau', A\tau\right\rangle   +\left\langle \varepsilon A^2 \tau, A\tau\right\rangle $$ $$ = \lambda \left( \beta_1 (t,x,v,\tau) \tau + \gamma(t,x,v,\tau) v, A\tau\right) . \ee
Note that we can use parentheses instead of brackets in the right-hand side due to the equality $$\langle w_1, w_2 \rangle= ( w_1, w_2 ),\ w_1\in L_2, w_2 \in H^1.$$

But \be\left\langle  \tau', A\tau\right\rangle=\left\langle  \tau',\tau\right\rangle_1=\frac 1 2\frac {d}{dt}\|\tau\|^2_1\ee (e.g. by \cite[Lemma  III.1.2]{tem}).
Thus,
\be \frac 1 2\frac {d}{dt}\|\tau\|^2_1   +\varepsilon(A \tau, A\tau)_1 $$ $$ = \lambda \left(\beta_1 (t,x,v,\tau) \tau + \gamma(t,x,v,\tau) v, \tau\right )+\lambda \left\langle \nabla[\beta_1 (t,x,v,\tau) \tau + \gamma(t,x,v,\tau) v], \nabla \tau\right \rangle $$ $$=\lambda \left(\beta_1 (t,x,v,\tau)\tau, \tau\right ) + \lambda\left (\gamma(t,x,v,\tau) v, \tau\right )$$ $$+\lambda \left ( \beta(t,x,v,\tau) \nabla v + \mu(t,x,v,\tau) \nabla \tau +g(t,x,v,\tau), \nabla \tau\right ) .\ee
Denote $\bar{v}(t)=e^{-kt}v(t)$, $\bar{\tau}(t)=e^{-kt}\tau(t)$, where $k>0$ will be defined below. Then
\be \frac 1 2\frac {d}{dt}\|e^{kt}\bar{\tau}\|^2_1  +e^{2kt}\varepsilon(A \bar\tau, A\bar\tau)_1 $$ $$ =\lambda \Big( \beta_1(t,x,e^{kt}\bar v(t),e^{kt}\bar\tau(t))  \bar \tau e^{kt},\bar \tau(t) e^{kt}\Big) + \lambda \Big( \gamma(t,x,e^{kt}\bar v(t),e^{kt}\bar\tau(t))  \bar ve^{kt},\bar \tau(t)e^{kt}\Big)$$ $$+\lambda \Big( \beta(t,x,e^{kt}\bar v(t),e^{kt}\bar\tau(t)) \nabla \bar ve^{kt} + \mu(t,x,e^{kt}\bar v(t),e^{kt}\bar \tau(t)) \nabla \bar\tau e^{kt} $$ $$+g(t,x,e^{kt}\bar v(t),e^{kt}\bar\tau(t)), \nabla\bar \tau(t)e^{kt}\Big) .\ee
Denote now $$\beta_k(t,x,\bar v(t),\bar\tau(t))=\beta(t,x,e^{kt}\bar v(t),e^{kt}\bar\tau(t)),$$
$$\mu_k(t,x,\bar v(t),\bar\tau(t))=\mu(t,x,e^{kt}\bar v(t),e^{kt}\bar\tau(t)),$$ $$g_k(t,x,\bar v(t),\bar\tau(t))=e^{-kt}g(t,x,e^{kt}\bar v(t),e^{kt}\bar\tau(t)),$$ $$\beta_{1k}(t,x,\bar v(t),\bar\tau(t))=\beta_1(t,x,e^{kt}\bar v(t),e^{kt}\bar\tau(t)),$$
$$\gamma_k(t,x,\bar v(t),\bar\tau(t))=\gamma(t,x,e^{kt}\bar v(t),e^{kt}\bar\tau(t)).$$

Thus,
\be \frac 1 2\frac {d}{dt}\|\bar{\tau}\|^2_1 + k \|\bar{\tau}\|^2_1    + \varepsilon(\bar\tau, \bar\tau)_X $$ $$ =\lambda \Big(\beta_{1k}(t,x,\bar v(t),\bar\tau(t))  \bar \tau(t), \bar \tau(t)\Big) + \lambda\Big(\gamma_k(t,x,\bar v(t),\bar \tau(t))  \bar v(t), \bar \tau(t)\Big) $$ $$+\lambda \Big(\beta_k(t,x,\bar v(t),\bar\tau(t)) \nabla \bar v + \mu_k(t,x,\bar v(t),\bar \tau(t)) \nabla \bar\tau $$ $$+g_k(t,x,\bar v(t),\bar\tau(t)), \nabla\bar \tau(t)\Big) .\ee
Integration from $0$ to $t\in[0,T]$ yields
\be   \frac 1 2\|\bar\tau(t)\|^2_1  +k\int\limits_0^t \|\bar{\tau}(s)\|^2_1\,ds+ \varepsilon\int\limits_0^t\|\bar\tau(s)\|^2_X\,ds  $$ $$ = \frac 1 2 \|\varsigma_0\|^2_1+ \lambda \int\limits_0^t \Big( \beta_k(s,x,\bar v(s),\bar \tau(s)) \nabla \bar v(s),  \bar \nabla\tau(s)\Big)+ \Big(\gamma_k(s,x,\bar v(s),\bar \tau(s))  \bar v(s), \bar \tau(s)\Big)$$ $$ +  \Big(\mu_k(s,x,\bar v(s),\bar\tau(s)) \nabla \bar\tau(s), \bar \nabla\tau(s)\Big)+\Big(\beta_{1k}(s,x,\bar v(s),\bar\tau(s))  \bar \tau(s), \bar \tau(s)\Big) $$ $$+ \Big(g_k(s,x,\bar v(s),\bar \tau(s)), \nabla\bar\tau(s)\Big)\,ds.  \ee
Applying the Cauchy-Buniakowski inequality, Cauchy's inequality $ab\leq c a^2+ \frac 1 {4c} b^2$, (3.15) and (3.16) we obtain
\be   \frac 1 2\|\bar\tau(t)\|^2_1  +k\int\limits_0^t \|\bar{\tau}(s)\|^2_1\,ds+ \varepsilon\int\limits_0^t\|\bar\tau(s)\|^2_X\,ds  $$ $$ \leq \frac 1 2 \|\varsigma_0\|^2_1+\frac{\lambda K^2_\beta}4\int\limits_0^t \|\bar{v}(s)\|^2_1\,ds +\lambda \int\limits_0^t \|\bar{\tau}(s)\|^2_1\,ds+\lambda K_\mu\int\limits_0^t \|\bar{\tau}(s)\|^2_1\,ds$$ $$+\frac {\lambda} 4 \int\limits_0^t \|g_k(s,\cdot,\bar v(s,\cdot),\bar \tau(s,\cdot))\|^2\,ds+ \lambda\int\limits_0^t \|\bar{\tau}(s)\|^2_1\,ds.  \ee
Note that $$\int\limits_0^t \|g_k(s,\cdot,\bar v(s,\cdot),\bar \tau(s,\cdot))\|^2\,ds\leq  \int\limits_0^t \|K_g[|\bar v(s,\cdot)|+|\bar \tau(s,\cdot)|]+\widetilde g(s,\cdot)\|^2\,ds$$ $$\leq 3 K_g^2\int\limits_0^t \|\bar v(s,\cdot)\|^2\,ds + 3 K_g^2\int\limits_0^t \|\bar \tau(s,\cdot)\|^2\,ds + 3\int\limits_0^t \|\widetilde g(s,\cdot)\|^2\,ds$$ $$\leq 3 K_g^2 \int\limits_0^t \|\bar v(s)\|^2_1\,ds + 3 K_g^2  \int\limits_0^t \|\bar \tau(s)\|_1^2\,ds + 3\|\widetilde g\|^2_{L_2((0,T)\times \Omega)}.$$
Hence,\be   \frac 1 2\|\bar\tau(t)\|^2_1  +(k-2-K_\mu-\frac 3 4 K_g^2 )\int\limits_0^t \|\bar{\tau}(s)\|^2_1\,ds+ \varepsilon\int\limits_0^t\|\bar\tau(s)\|^2_X\,ds  $$ $$ \leq \frac 1 2 \|\varsigma_0\|^2_1+\lambda(\frac {K^2_\beta} 4+\frac 3 4 K_g^2 )\int\limits_0^t \|\bar{v}(s)\|^2_1\,ds +\frac 34\|\widetilde g\|^2_{L_2((0,T)\times \Omega)}.  \ee

Take $k\geq 4+2K_\mu+\frac 3 2 K_g^2$.

In particular, (4.14) implies
\be   \int\limits_0^t \|\bar{\tau}(s)\|^2_1\,ds\leq \frac {C} k(1+\lambda\int\limits_0^t \|\bar{v}(s)\|^2_1\,ds).  \ee

Now, take the "bra-ket" of (4.1) (as
elements of $H^{-2}_N(\Omega)$) and $v(t)\in H^2(\Omega)$ at a.a.
$t\in[0,T]$: \be \left\langle  v', v\right\rangle   + \left\langle \varepsilon A_2 v, v\right\rangle $$ $$ = \lambda \left\langle div_N[D(t,x,v,\tau) \nabla v + E(t,x,v,\tau) \nabla \tau +f(t,x,v,\tau)]+\psi, v\right\rangle . \ee
Again, by \cite[Lemma  III.1.2]{tem}, \be \left\langle  v',v\right\rangle=\frac 1 2\frac {d}{dt}\|v\|^2.\ee
Thus,
\be \frac 1 2\frac {d}{dt}\|v\|^2   + \varepsilon(v, v)_2 $$ $$ =-\lambda (D(t,x,v,\tau) \nabla v + E(t,x,v,\tau) \nabla \tau +f(t,x,v,\tau), \nabla v) +\lambda\left\langle\psi, v\right\rangle.\ee
Then
\be \frac 1 2\frac {d}{dt}\|e^{kt}\bar{v}\|^2   + e^{2kt}\varepsilon( \bar v,  \bar v )_2 $$ $$ =-\lambda \Big( D(t,x,e^{kt}\bar v(t),e^{kt}\bar\tau(t)) \nabla \bar ve^{kt} + E(t,x,e^{kt}\bar v(t),e^{kt}\bar \tau(t)) \nabla \bar\tau e^{kt} $$ $$+f(t,x,e^{kt}\bar v(t),e^{kt}\bar\tau(t)), \nabla\bar v(t)e^{kt}\Big) +\lambda\left\langle\psi(t),\bar v(t)e^{kt}\right\rangle.\ee
Denote now $$D_k(t,x,\bar v(t),\bar\tau(t))=D(t,x,e^{kt}\bar v(t),e^{kt}\bar\tau(t)),$$
$$E_k(t,x,\bar v(t),\bar\tau(t))=E(t,x,e^{kt}\bar v(t),e^{kt}\bar\tau(t)),$$ $$f_k(t,x,\bar v(t),\bar\tau(t))=e^{-kt}f(t,x,e^{kt}\bar v(t),e^{kt}\bar\tau(t)).$$

Thus, \be \frac 1 2\frac {d}{dt}\|\bar{v}\|^2 + k \|\bar{v}\|^2 +
\varepsilon(\bar v, \bar v )_2 $$ $$ =-\lambda \Big(D_k(t,x,\bar
v(t),\bar\tau(t)) \nabla \bar v + E_k(t,x,\bar v(t),\bar \tau(t))
\nabla \bar\tau $$ $$+f_k(t,x,\bar v(t),\bar\tau(t)), \nabla\bar
v(t)\Big) + e^{-kt}\lambda\left\langle\psi(t),\bar
v(t)\right\rangle.\ee Therefore \be   \frac 1 2\|\bar v(t)\|^2
+k\int\limits_0^t \|\bar{v}(s)\|^2\,ds+
\varepsilon\int\limits_0^t\|\bar v(s)\|^2_2\,ds  $$ $$ = \frac 1 2
\|u_0\|^2- \lambda \int\limits_0^t\Big( D_k(s,x,\bar v(s),\bar
\tau(s)) \nabla \bar v(s) + E_k(s,x,\bar v(s),\bar\tau(s)) \nabla
\bar\tau(s) $$ $$+f_k(s,x,\bar v(s),\bar \tau(s)), \nabla\bar
v(s)\Big)- e^{-ks}\left\langle\psi(s),\bar v(s)\right\rangle\,ds.
\ee Using Cauchy's inequality, (3.14) and (3.19), we get \be \frac
1 2\|\bar v(t)\|^2  +k\int\limits_0^t \|\bar{v}(s)\|^2\,ds+
\varepsilon\int\limits_0^t\|\bar v(s)\|^2_2\,ds  + \lambda d
\int\limits_0^t (\nabla \bar v(s), \nabla\bar v(s))\,ds$$ $$ \leq
\frac 1 2 \|u_0\|^2+\frac {\lambda K^2_E}{d}\int\limits_0^t
\|\nabla\bar{\tau}(s)\|^2\,ds + \frac{\lambda d}{4}\int\limits_0^t
\|\nabla\bar{v}(s)\|^2\,ds$$ $$+\frac {\lambda} {d}
\int\limits_0^t \|f_k(s,\cdot,\bar v(s,\cdot),\bar
\tau(s,\cdot))\|^2\,ds+  \frac{\lambda d}{4}\int\limits_0^t
\|\nabla\bar{v}(s)\|^2\,ds. $$ $$+\frac {\lambda} {d}
\int\limits_0^t \|\psi(s)\|^2_{-1}\,ds+  \frac{\lambda
d}{4}\int\limits_0^t \|\bar{v}(s)\|^2_1\,ds.  \ee As for $g_k$
above, we have $$\int\limits_0^t \|f_k(s,\cdot,\bar
v(s,\cdot),\bar \tau(s,\cdot))\|^2\,ds$$ $$  \leq 3
K_f^2\int\limits_0^t \|\bar v(s)\|^2\,ds + 3 K_f^2 \int\limits_0^t
\|\bar \tau(s)\|_1^2\,ds + 3\|\widetilde f\|^2_{L_2((0,T)\times
\Omega)}.$$ Hence, from (4.22) and (4.15), \be   \frac 1 2\|\bar
v(t)\|^2  +(k-\frac {3 K_f^2}{d}-\frac d 4)\int\limits_0^t
\|\bar{v}(s)\|^2\,ds+ \varepsilon\int\limits_0^t\|\bar
v(s)\|^2_2\,ds  + \frac{\lambda d}{4} \int\limits_0^t \|\nabla\bar
v(s)\|^2\,ds$$ $$\leq  \frac 1 2 \|u_0\|^2+(\frac {K^2_E}{d}+
\frac 3{d}K_f^2 )\int\limits_0^t \|\bar{\tau}(s)\|^2_1\,ds +\frac
3 {d} \|\widetilde f\|^2_{L_2((0,T)\times \Omega)}+\frac 1 {d}
\|\psi\|^2_{L_2(0,T;H^{-1}_N (\Omega))}.  $$ $$ \leq \frac {C_0}
k(1+\lambda\int\limits_0^t \|\bar{v}(s)\|^2_1\,ds)+C$$ $$=\frac
{C_0} k(1+\lambda\int\limits_0^t \|\bar{v}(s)\|^2\,ds +
\lambda\int\limits_0^t \|\nabla\bar{v}(s)\|^2\,ds)+C.\ee

Take $k\geq \frac {3 K_f^2}{d}+\frac {8C_0}d+\frac d 4+\frac {C_0}k+1$.  Then (4.23) yields $$\int\limits_0^t \|\bar{v}(s)\|^2\,ds+\frac{\lambda d} 8\int\limits_0^t \|\nabla\bar{v}(s)\|^2\,ds \leq C$$
(now $C$ may depend on $k$), so $$ {\lambda }\int\limits_0^t \|\bar{v}(s)\|^2_1\,ds \leq C.$$  Thus, the right-hand members of inequalities (4.14) and (4.23) are bounded, and we arrive at
\be
\varepsilon\|\bar v\|^2_{L_2(0,T;H^2(\Omega))}+ \varepsilon\|\bar\tau\|^2_{L_2(0,T;X_N)}+$$ $$ \|\bar v\|^2_{L_\infty(0,T;L_2(\Omega))} + \lambda\|\bar v\|^2_{L_2(0,T;H^1(\Omega))} +
\|\bar \tau\|^2_{L_\infty(0,T;H^1(\Omega))}\leq C. \ee
Since $e^{kt} \leq e^{kT}$ for $t\in [0,T]$, this implies (4.6).
$ \Box $

\begin{lemma} Let $(v,\tau)$ be a weak solution to problem
(4.1)-(4.4). Then there is the following bound of the time derivatives: \be
 \|v'\|_{L_2(0,T;H^{-2}_N(\Omega))} +
\|\tau'\|_{L_2(0,T;H^{-1}_N(\Omega))}\leq C(1+\sqrt{\varepsilon}) \ee where $C$ is
independent of $\lambda$ and $\varepsilon$. \end{lemma}
\textbf{Proof.} Really, since $H^{-1}_N(\Omega)\subset H^{-2}_N(\Omega)$ continuously, (4.1) and (4.6) imply $$\|v'\|_{L_2(0,T;H^{-2}_N(\Omega))}\leq \varepsilon \|A_2 v\|_{L_2(0,T;H_N^{-2}(\Omega))}+$$ $$ \lambda  \|div_N[D(t,x,v,\tau) \nabla v + E(t,x,v,\tau) \nabla \tau +f(t,x,v,\tau)]\|_{L_2(0,T;H^{-1}_N(\Omega))}+\lambda\|\psi\|_{L_2(0,T;H^{-1}_N(\Omega))}$$
$$\leq \sqrt\varepsilon \sqrt\varepsilon\| v\|_{L_2(0,T;H^{2}(\Omega))}+$$ $$\lambda  \|D(t,x,v,\tau) \nabla v + E(t,x,v,\tau) \nabla \tau +f(t,x,v,\tau)\|_{L_2(0,T;L_2(\Omega))}+\|\psi\|_{L_2(0,T;H^{-1}_N(\Omega))}$$ $$\leq C\sqrt\varepsilon + K_D\lambda \| v\|_{L_2(0,T;H^1(\Omega))} + K_E\|\tau\|_{L_2(0,T;H^1(\Omega))} +\|f(t,x,v,\tau)\|_{L_2(0,T;L_2(\Omega))}+C$$ $$\leq C\sqrt\varepsilon + K_D\sqrt\lambda \| v\|_{L_2(0,T;H^1(\Omega))} + K_E\|\tau\|_{L_2(0,T;H^1(\Omega))} $$ $$+K_f\|v\|_{L_2(0,T;L_2(\Omega))}+K_f\|\tau\|_{L_2(0,T;L_2(\Omega))}+\|\widetilde f\|_{L_2((0,T)\times\Omega)}+C$$ $$\leq C\sqrt\varepsilon + C [\sqrt\lambda \| v\|_{L_2(0,T;H^1(\Omega))} +\|\tau\|_{L_\infty(0,T;H^1(\Omega))} $$ $$+\|v\|_{L_\infty(0,T;L_2(\Omega))}+\|\tau\|_{L_\infty(0,T;H^1(\Omega))}+1]\leq C(1+\sqrt\varepsilon).$$

Similarly, since $L_2(\Omega)\subset H^{-1}_N(\Omega)$ continuously, (4.2) and (4.6) yield $$\|\tau'\|_{L_2(0,T;H^{-1}_N(\Omega))}\leq \varepsilon \|A^2 \tau\|_{L_2(0,T;H^{-1}_N(\Omega))}+$$ $$ \lambda \|\beta_1 (t,x,v,\tau) \tau + \gamma(t,x,v,\tau) v\|_{L_2(0,T;L_2(\Omega))}$$
$$\leq \sqrt\varepsilon \sqrt\varepsilon\| \tau\|_{L_2(0,T;X_N)}+$$ $$K_\mu\| \tau\|_{L_2(0,T;L_2(\Omega))} + K_\beta\| v\|_{L_2(0,T;L_2(\Omega))}\leq C(1+\sqrt\varepsilon).$$$ \Box $

\begin{lemma} Given $u_0 \in L_2(\Omega)$, $\varsigma_0\in H^1(\Omega)$, there exists a weak solution to problem (4.1)-(4.4) in class (4.5). \end{lemma}

 \textbf{Proof.} Let us introduce auxiliary operators by the following formulas:
$$Q_1:W_2\times W_1\to L_2 (0, T; H^{-2}_N(\Omega)), $$
$$Q_1 (v, \tau) = div_N[D(\cdot,\cdot,v,\tau) \nabla v],$$
$$Q_2:W_2\times W_1\to L_2 (0, T; H^{-2}_N(\Omega)), $$
$$Q_2 (v, \tau) = div_N[E(\cdot,\cdot,v,\tau) \nabla \tau],$$
$$Q_3:W_2\times W_1\to L_2 (0, T; H^{-2}_N(\Omega)), $$
$$Q_3 (v, \tau) = div_N[f(\cdot,\cdot,v,\tau)]+\psi,$$
$$Q_4:W_2\times W_1\to L_2 (0, T; H^{-1}_N(\Omega)), $$
$$Q_4 (v, \tau) = \gamma(\cdot,\cdot,v,\tau) v,$$
$$Q_5:W_2\times W_1\to L_2 (0, T; H^{-1}_N(\Omega)), $$
$$Q_5 (v, \tau) = \beta_1(\cdot,\cdot,v,\tau) \tau,$$

$$Q:W_2\times W_1\to L_2 (0, T; H^{-2}_N(\Omega))\times L_2 (0, T; H ^ {-1}_N(\Omega))\times L_2(\Omega)\times H^1(\Omega), $$
$$Q (v, \tau) = (-Q_1 (v,\tau) -Q_2(v, \tau) -Q_3(v,\tau), -Q_4 (v,\tau) -Q_5(v, \tau), 0,0), $$

$$\tilde{A}_1:W_1\to L_2 (0, T; H^{-1}_N(\Omega))\times  H^1(\Omega),$$ $$
\tilde{A}_1 (u)= (u'+\varepsilon A^2 u, u | _ {t=0})
,
$$
$$\tilde{A}_2:W_2\to L_2 (0, T; H^{-2}_N(\Omega))\times  L_2(\Omega),$$ $$
\tilde{A}_2 (u)= (u'+\varepsilon A_2
u, u | _ {t=0}),
$$
$$\tilde{A}:W_2\times W_1\to L_2 (0, T; H^{-2}_N(\Omega))\times L_2 (0, T; H ^ {-1}_N(\Omega))\times L_2(\Omega)\times H^1(\Omega),$$ $$
\tilde{A} (v,\tau)= (v'+\varepsilon A_2
v, \tau'+\varepsilon A^2
\tau , v | _ {t=0}, \tau | _ {t=0}).$$

 Then the weak statement of problem (4.1) - (4.4) is equivalent to the operator
 equation
\begin{equation}\tilde {A} (v, \tau) + \lambda Q (v, \tau) = (0, 0, u_0, \varsigma_0). \end{equation}

Let us briefly explain the idea of the proof. We are going to show that the operator $\tilde {A}$ is invertible. This yields the solvability of  equation (4.26) for $\lambda=0$. On the other hand, $Q$ turns out to be a compact operator. Then we can rewrite (4.26) in a form suitable for application of the Leray-Schauder degree theory, which will imply the existence of solutions for all $\lambda\in [0,1]$.

We recall that a non-linear operator $K:X_1\to X_2$ ($X_1$ and $X_2$ are Banach spaces) is called compact if it is continuous and the image of any bounded set in $X_1$ is relatively compact in $X_2$.
In particular, if $X_1$ is reflexive, and, for any sequence $x_m\to x_*$ which converges in $X_1$ in the weak sense, one has $K(x_m)\to K(x_*)$ strongly in $X_2$, then $K$ is compact (since any bounded subset of $X_1$
is relatively compact in the weak topology).

For some $q>2$, the embeddings $W_1\subset L_q (0, T;{W}{}_q^1(\Omega)), $ $W_2\subset
L_q (0, T;{W}{}_q^1(\Omega)) $ are compact. Really, we have $W_1\subset C([0,T];H^1(\Omega)),\ W_2\subset
C([0,T];L_2(\Omega))$ continuously.
Note that (by
the Rellich-Kondrashov theorem) $H^2 \subset L_2$ compactly. Furthermore,  $H^1 \subset H^{-1}_N$ compactly, so the adjoint embedding $X_N\subset H^1$ is also compact.
Then, by \cite[Corollary 6]{sim1}, $W_1\subset L_p(0,T;H^1(\Omega)),\ W_2\subset
L_p(0,T;L_2(\Omega))$ compactly for every $p<\infty$. Let $p_1=\frac {2n} {n-1}$  and $p_0= \frac {2n} {n-2}$, cf. (2.7). Then $\frac 2 {p_1}=\frac 1 2 + \frac 1 {p_0}$.
For $u\in X_N$, we have
$$\|u\|^2_{W^1_{p_1}}\leq C(\|u\|^2_{L_{p_1}}+\|\nabla u\|^2_{L_{p_1}}).$$ The second term is $$\||\nabla u|^2\|_{L_{p_1/2}}\leq \|\nabla u\|_{L_{2}}\|\nabla u\|_{L_{p_0}}\leq C\|u\|_1\|u\|_X.$$ The first term can be estimated similarly.
If $q_1>2$ is such that $\frac 2 {q_1}=\frac 1 2 + \frac 1 {p}$ with some $p$ large enough, then, by
 \cite[Lemma 11]{sim1}, $W_1\subset L_{q_1}(0,T;W^1_{p_1})$ compactly. Now, for $u\in H^2$, we have
$$\|u\|^4_{W^1_{p_1}}\leq C(\|u\|^4_{L_{p_1}}+\|\nabla u\|^4_{L_{p_1}}),$$  and the second term is $$\||\nabla u|^2\|^2_{L_{p_1/2}}\leq \|\nabla u\|^2_{L_{2}}\|\nabla u\|^2_{L_{p_0}}$$ $$\leq C\|u\|_1^2\|u\|^2_{2}\leq C\|u\|\|u\|^3_2.$$ The last inequality follows from \cite[Theorem 4.17]{adms}. The first term can be estimated in the same way.
If $q_2>2$ is such that $\frac 4 {q_2}=\frac 3 2 + \frac 1 {p}$ with some $p$ large enough, then, by
 \cite[Lemma 11]{sim1}, $W_2\subset L_{q_2}(0,T;W^1_{p_1})$ compactly.

Let us show that the operators $Q_1,\dots,
Q_5$ are compact. Let $v_m\to v_*$ weakly in $W_2$, $\tau_m\to \tau_*$ weakly in $W_1$.
Then $v_m\to v_*$, $\tau_m\to \tau_*$ strongly in $L_q (0, T; W^1_q(\Omega))$ and in $L_q (0, T; L_q(\Omega))$, and $\nabla v_m\to \nabla v_*$, $\nabla \tau_m\to \nabla\tau_*$ strongly in $L_q (0, T; L_q(\Omega)^n)$.

By Krasnoselskii's theorem \cite{kras,skry} on
continuity of Nemytskii operators we have
$$D(\cdot,\cdot,v_m,\tau_m) \to D(\cdot,\cdot,v_*,\tau_*),$$
$$E(\cdot,\cdot,v_m,\tau_m) \to E(\cdot,\cdot,v_*,\tau_*),$$
strongly in $L_p ((0, T)\times \Omega)^{n\times n}$,
$$\beta_1(\cdot,\cdot,v_m,\tau_m) \to \beta_1(\cdot,\cdot,v_*,\tau_*),$$
$$\gamma(\cdot,\cdot,v_m,\tau_m) \to \gamma(\cdot,\cdot,v_*,\tau_*),$$
strongly in $L_p ((0, T)\times \Omega)$ for all $p<\infty$, and
$$f(\cdot,\cdot,v_m,\tau_m) \to f(\cdot,\cdot,v_*,\tau_*),$$
strongly in $L_2 ((0, T)\times \Omega)^n$.

Clearly, if a sequence of functions $y_m$ converges in $L_q ((0, T)\times \Omega)$, and another sequence $z_m$ converges in $L_p ((0, T)\times \Omega)$, $\frac 1 p + \frac 1 q =\frac 12$, then their pointwise products $y_m z_m$ tend to the product of their limits in $L_2 ((0, T)\times \Omega)$.

Hence, $D(\cdot,\cdot,v_m,\tau_m) \nabla  v_m\to D(\cdot,\cdot,v_*,\tau_*)\nabla v_*$ in $L_2 (0, T; L_2(\Omega)^n)$. Therefore $$Q_1(v_m,\tau_m)\to Q_1(v_*,\tau_*)$$ in $L_2 (0, T; H^{-1}(\Omega))$ (and all the more in $L_2 (0, T; H^{-2}(\Omega))$). Similarly, $$Q_i(v_m,\tau_m)\to Q_i(v_*,\tau_*),\, i=2,\dots,5,$$ in $L_2 (0, T; H^{-1}(\Omega))$.

Hence, the operator $Q $ is also compact.

Note that $$\left\langle A_2 u, u\right\rangle =\|u\|_2^2,$$ for $u\in H^2(\Omega)$, and
$$\left\langle A^2 u, u\right\rangle_1=\left\langle A^2 u, A u\right\rangle=(A u, A u)_1=\|u\|_X^2$$ for $u\in X$.
Therefore the operators $
\tilde {A}_1 $ and $\tilde {A}_2$ are invertible (e.g. by Theorem 1.1 from \cite{ggz},
Chapter VI, or Lemma 3.1.3 from \cite{gruy}). Hence, $
\tilde {A}$ is also (continuously) invertible.

Rewrite equation (4.26) as

\begin{equation} (u, \tau) + \lambda\tilde {A} ^ {-1} Q (u, \tau) = \tilde {A} ^ {-1} (0, 0, u_0,
\varsigma_0) .\end{equation} A priori bounds from Lemmas 4.1 and 4.2 imply that equation (4.27) has no
solutions on the boundary of a sufficiently large ball $B $ in
$W_2\times W_1, $ independent of $ \lambda. $ Without loss of
generality $a_0 =\tilde {A} ^ {-1} (0, 0, u_0, \varsigma_0) $ belongs to
this ball. Then we can consider the Leray - Schauder degree (see
e.g. \cite{ll}) of the map $I + \lambda \tilde {A} ^ {-1} Q $ ($I$
is the identity map) on the ball $B $ with respect to the point
$a_0 $,
$$deg _ {LS} (I + \lambda \tilde {A} ^ {-1} Q, B, a_0). $$  By the homotopic invariance property of
the degree we have
$$ deg _ {LS} (I + \lambda \tilde {A} ^ {-1} Q, B, a_0) = deg _ {LS} (I, B,
a_0) =1\neq 0 .$$ Thus, equation (4.27) (and, therefore, problem
(4.1) - (4.4)) has a solution in the ball $B $ for every $
\lambda. $ $ \Box $

\textbf{Proof of Theorem 3.1.} Take a decreasing sequence of
positive numbers $\varepsilon_m\to 0$.  By Lemma 4.3, there
is a pair $ (v_m, \tau_m) $ which is a weak solution to problem
(4.1)-(4.4) with $\lambda=1$,
$\varepsilon=\varepsilon_m$.

Due to a priori estimate
(4.6), without loss of generality (passing to a subsequence if
necessary) one may assume that there exist
limits\\
 $u =\lim\limits _ {m\to
\infty} ^ {}  {v} _ {m} $, which is $*$-weak in
$L _ {\infty} (0, T; L_2(\Omega)) $ and weak in $L _ {2} (0, T; H^1(\Omega));$ \\
$ \varsigma =\lim\limits _ {m\to \infty} ^ {} {\tau} _ {m}
$, which is $ *$-weak in $L _ {\infty} (0, T; H^1(\Omega)) $ and weak in
$L _ {2} (0, T; H^1(\Omega))$.

Moreover, due to Lemma 4.2, without loss of generality one may
assume that ${v}_m '\to u '$ weakly in $L_2(0,T;H^{-2}_N)$, $
{\tau}_m ' \to \varsigma ' $ weakly in $L_2(0,T;H^{-1}_N)$. Then,
by  \cite[Corollary 4]{sim1}, $v_m\to u,$ $\tau_m \to \varsigma$
strongly in $C ([0, T]; H^{-1}_N)$. Therefore $u$ and $\varsigma$
satisfy (3.8).

Furthermore, by  \cite[Corollary 4]{sim1},
$v_m\to u,$
$\tau_m \to \varsigma$
strongly in $L_2 (0, T; L_2)$.

By Krasnoselskii's theorem \cite{kras,skry} we have again
$$D(\cdot,\cdot,v_m,\tau_m) \to D(\cdot,\cdot,u,\varsigma),$$
$$E(\cdot,\cdot,v_m,\tau_m) \to E(\cdot,\cdot,u,\varsigma),$$
strongly in $L_p ((0, T)\times \Omega)^{n\times n}$,
$$\beta_1(\cdot,\cdot,v_m,\tau_m) \to \beta_1(\cdot,\cdot,u,\varsigma),$$
$$\gamma(\cdot,\cdot,v_m,\tau_m) \to \gamma(\cdot,\cdot,u,\varsigma),$$
strongly in $L_p ((0, T)\times \Omega)$ for all $p<\infty$, and
$$f(\cdot,\cdot,v_m,\tau_m) \to f(\cdot,\cdot,u,\varsigma),$$
strongly in $L_2 ((0, T)\times \Omega)^n$.

Observe that if a sequence of functions $y_m$ converges weakly in $L_2 ((0, T)\times \Omega)$, and another sequence $z_m$ converges strongly in $L_p ((0, T)\times \Omega)$, then their pointwise products $y_m z_m$ converge weakly to the product of their limits in $L_q ((0, T)\times \Omega)$, $\frac 1 p + \frac 1 2=\frac 1 q $.

Therefore, $$D(\cdot,\cdot,v_m,\tau_m) \nabla  v_m\to D(\cdot,\cdot,u,\varsigma)\nabla u,$$
$$E(\cdot,\cdot,v_m,\tau_m) \nabla  \tau_m\to E(\cdot,\cdot,u,\varsigma)\nabla \varsigma,$$
weakly in $L_q (0, T; L_q(\Omega)^n)$,
$$\gamma(\cdot,\cdot,v_m,\tau_m) v_m\to \beta(\cdot,\cdot,u,\varsigma)u,$$
$$\beta_1(\cdot,\cdot,v_m,\tau_m) \tau_m\to \beta_1(\cdot,\cdot,u,\varsigma)\varsigma$$ weakly in $L_q (0, T; L_q(\Omega))$ for $1\leq q<2$. Therefore the right-hand members of (4.1) converge to the corresponding right-hand members of (3.21) weakly in $L_q (0, T; W^{-1}_{q,N}(\Omega))$.

Due to Lemma 4.1, $\varepsilon_m\|v_m\|_{L_2(0,T;H^2(\Omega))}=\sqrt{\varepsilon_m}\sqrt{\varepsilon_m}\|v_m\|_{L_2(0,T;H^2(\Omega))}\to 0.$ Hence, $\varepsilon_m\|A_2 v_m\|_{L_2(0,T;H^{-2}_N(\Omega))}\to 0.$ Similarly, $\varepsilon_m\|\tau_m\|_{L_2(0,T;X_N)}$ tends to zero, so\\ $\varepsilon_m\|A^2\tau_m\|_{L_2(0,T;H_N^{-1})}\to 0$.

W.l.o.g. we may assume, in addition, that $q\geq \frac {2n}
{n+2}$. Then, by Sobolev theorem,  $H^2(\Omega)\subset
W^1_{{q}/{q-1}}(\Omega)$, so $W^{-1}_{q}(\Omega)\subset
H^{-2}_N(\Omega)$. Passing to the limit as $m \to \infty$ in (4.1)
and (4.2) with $\lambda=1$, $\varepsilon=\varepsilon_m$, $v=v_m$,
$\tau=\tau_m$ in the space of distributions on $(0,T)$ with values
in  $H^{-2}_N(\Omega)$ (for (4.2) it is possible for $H^{-1}_N$ as
well), we conclude that the pair $(u,\varsigma)$ is a solution to
(3.5)--(3.8).

It remains to observe that the right-hand side (and, hence, the
left-hand side) of  (3.21) belongs to $L_2(0,T;H^{-1}_N)$, and,
due to (3.9), the ones of (3.6) belong to $L_2(0,T;H^{1})$. $ \Box
$

\section{Long-time behaviour}

Theorem 3.1 implies that solutions to (3.5)--(3.8) can be continued, step by step, onto the whole positive semi-axis:
\begin{cor} Given $u_0 \in L_2(\Omega)$, $\varsigma_0\in H^1(\Omega)$ and $\psi\in L_{2,loc}(0,\infty;H^{-1}_N(\Omega))$, there is a pair  \be u \in L_{2,loc}(0,\infty;H^1(\Omega))\bigcap H^1_{loc}(0,\infty;H^{-1}_N(\Omega)), \varsigma\in
H^1_{loc}(0,\infty;H^1(\Omega))\ee which satisfies (3.8), whereas (3.21) holds true in $H^{-1}_N(\Omega)$ a.e. on
$(0,\infty)$, and (3.6) holds a.e. in
$(0,\infty)\times \Omega$. \end{cor}

Below we keep assuming conditions i)-vi) of
Section 3, but we replace (3.17)--(3.19)
with stronger requirements, namely

vii)
\be |f (t, x, v, \tau) | \leq \widetilde f (t, x), |g (t, x, v, \tau) | \leq \widetilde g (t, x) \ee
with some known functions $ \widetilde f, \widetilde g \in L_2 ((0, \infty) \times\Omega) $, and

viii) there are\footnote{See \cite[Section 5]{iter} for a discussion whether this assumption is realistic.} positive numbers $ \Gamma $ and $ \Gamma_0 $ such that \be (D (\cdot) \xi, \xi) _ {\R^n} - (\mu (\cdot) \eta, \eta) _ {\R^n} + \left (\left [E (\cdot) \Gamma-\frac {\beta (\cdot)} {\Gamma} \right] \xi, \eta\right) _ {\R^n} \geq \Gamma_0 (| \xi | ^ 2 + |\eta | ^ 2) \ee
for any $ \xi, \eta\in \R^n $.

Consider any global weak solution $(u,\varsigma)$ existing by Corollary 5.1. Denote $\Psi=A^{-1}\psi$.
Then $\Psi\in L_{2,loc}(0,\infty;H^{1}(\Omega))$. Assume, in addition, that $\Psi\in L_{2}(0,\infty;H^{1}(\Omega))\cap L_{1}(0,\infty;L_2(\Omega))$.

\begin {lemma} The following estimate is valid:
\be
 \|u\|_{L_\infty (0, \infty;L_2( \Omega))} +
 \| \nabla u\| _{L_2 (0, \infty;L_2( \Omega))}  $$ $$+ \| \nabla\varsigma \|_{L_\infty (0, \infty;L_2( \Omega))} +
  \|\nabla\varsigma\|_{L_2 (0, \infty;L_2( \Omega))}
\leq C.\ee \end {lemma}

\textbf {Proof}. The condition (5.3) can be rewritten as \be (D
(\cdot) \Gamma^2\xi, \xi) _ {\R^n} - (\mu (\cdot) \eta, \eta) _
{\R^n} + \left (E (\cdot) \Gamma^2\xi, \eta\right) _ {\R^n} -
({\beta (\cdot)} \xi, \eta) _ {\R^n} $$ $$\geq \Gamma_0 (\Gamma^2
|\xi | ^ 2 + |\eta | ^ 2) \ee for any $ \xi, \eta\in \R^n $ (just
substitute  $ \Gamma\xi $ for $ \xi  $ in (5.3)).

Take the "bra-ket" of $-div_N \nabla \varsigma(t)\in
H^{-1}_N(\Omega)$ and the terms of (3.6)  (as elements of
$H^{1}(\Omega)$)  at a.a. $t\in[0,T]$: \be \left(\nabla
\varsigma', \nabla\varsigma\right)
 = \left(\nabla( \beta_1 (t,x,u,\varsigma) \varsigma +
\gamma(t,x,u,\varsigma) u), \nabla\varsigma\right) . \ee

Thus, \be \frac 1 2\frac {d}{dt}\|\nabla \varsigma\|^2   = \left (
\beta(t,x,u,\varsigma) \nabla u + \mu(t,x,u,\varsigma) \nabla
\varsigma +g(t,x,u,\varsigma), \nabla \varsigma\right ) .\ee

Take the "bra-ket" of (3.21) (as elements of $H^{-1}_N(\Omega)$)
and $u(t)\in H^1(\Omega)$ at a.a. $t\in[0,T]$, arriving at (cf.
the proof of Lemma 4.1) \be \frac 1 2\frac {d}{dt}\|u\|^2   =-
(D(t,x,u,\varsigma) \nabla u + E(t,x,u,\varsigma) \nabla \varsigma
+f(t,x,u,\varsigma), \nabla u) +\left\langle\psi,
u\right\rangle.\ee  Multiply it by $ \Gamma^2 $ and add this with
(5.7):
 \be \frac {\Gamma^2} 2\frac {d} {dt} \|u \| ^ 2 + \frac 1 2\frac {d} {dt} \| \nabla \varsigma \| ^ 2 $$ $$ =- (D (t, x, u, \varsigma) \Gamma^2 \nabla u + E (t, x, u, \varsigma) \Gamma^2\nabla \varsigma, \nabla u) $$ $$ + \Big (\beta (t, x, u, \varsigma) \nabla u + \mu (t, x, u, \varsigma) \nabla \varsigma, \nabla \varsigma\Big) $$ $$ + \Big (g (t, x, u, \varsigma), \nabla \varsigma\Big)-\Gamma^2 (f (t, x, u, \varsigma), \nabla u) +\Gamma^2\langle\psi, u\rangle.\ee

Using (5.5), we conclude that

\be \frac {\Gamma^2} 2\frac {d} {dt} \|u \| ^ 2 +   \frac 1 2\frac
{d} {dt} \| \nabla\varsigma \| ^ 2 + \Gamma_0 (\Gamma^2 \|\nabla u
\| ^ 2 + \|\nabla \varsigma \| ^ 2) $$ $$ \leq \Big (g (t, x, u,
\varsigma), \nabla \varsigma\Big)-\Gamma^2 (f (t, x, u,
\varsigma), \nabla u) +\Gamma^2(\Psi, u)_1 .\ee

Integrating along the interval $ (0, t) $, $t>0$, we get \be \frac
{\Gamma^2} 2 \|u(t) \| ^ 2 +   \frac 1 2\| \nabla\varsigma(t) \| ^
2 +
  \Gamma_0\Gamma^2\int\limits_0^t \| \nabla u (s) \| ^2 \, ds + \Gamma_0\int\limits_0^t \|\nabla\varsigma (s) \| ^2 \, ds $$ $$ \leq  \frac {\Gamma^2} 2 \|u_0 \| ^ 2 +   \frac 1 2\| \nabla\varsigma_0 \| ^ 2+ \int\limits_0^t\Big (g (s, x, u(s), \varsigma (s)), \nabla\varsigma (s) \Big) \, ds $$ $$+\int\limits_0^t\Gamma^2(\Psi(s), u(s))  \, ds$$ $$+\int\limits_0^t\left[\Gamma^2(\nabla\Psi(s), \nabla u(s))- \Gamma^2\Big (f (s, x, u (s), \varsigma (s)), \nabla u (s) \Big)\right] \, ds. \ee

Applying the Cauchy-Buniakowski inequality, Cauchy's inequality
and (5.2), we observe that \be\big | \int\limits_0^t\Big
(\nabla\Psi(s)+f (s, x, u (s), \varsigma (s)), \nabla u (s) \Big)
\, ds\big | $$ $$\leq [\|\nabla\Psi\|_ {L_2 ((0, \infty) \times
\Omega)}+\| \widetilde f \| _ {L_2 ((0, \infty) \times \Omega)}]
\| \nabla u \| _ {L_2 ((0, t) \times \Omega)} $$ $$\leq \frac 1
{2\Gamma_0} [\|\nabla\Psi\|_ {L_2 (0, \infty;L_2( \Omega))}+\|
\widetilde f \| _ {L_2 ((0, \infty) \times \Omega)}]^2 + \frac
{\Gamma_0} 2 \|\nabla u \| ^ 2 _ {L_2 (0, t;L_2( \Omega))} .\ee

Similarly, \be\big | \int\limits_0^t\Big (g (s, x, u (s),
\varsigma (s)), \nabla \varsigma (s) \Big) \, ds\big | \leq \frac
1 {2\Gamma_0} \| \widetilde g \| ^ 2 _ {L_2 ((0, \infty) \times
\Omega)} + \frac {\Gamma_0} 2 \|\nabla \varsigma \| ^ 2 _ {L_2 (0, t;L_2( \Omega))}. \ee

And, obviously, \be \big |\int\limits_0^t\Gamma^2(\Psi(s), u(s))  \, ds \big | \leq \Gamma^2\|\Psi\|_ {L_1 (0, \infty;L_2( \Omega))} \|u\|_{L_\infty (0, t;L_2( \Omega))}.\ee

Inequalities (5.11)--(5.14) yield \be \frac {\Gamma^2} 2 \|u(t) \|
^ 2 +   \frac 1 2\| \nabla\varsigma(t) \| ^ 2 +
  \frac{\Gamma_0\Gamma^2} 2 \int\limits_0^t \| \nabla u (s) \| ^2 \, ds + \frac{\Gamma_0} 2\int\limits_0^t \|\nabla\varsigma (s) \| ^2 \, ds $$ $$ \leq   C+C\|u\|_{L_\infty (0, t;L_2( \Omega))}, \ee where $C$ is independent of $t$,
thus we have the same inequality between the essential supremums
of both members on $(0,t)$, in particular, $$ \frac {\Gamma^2} 2
\|u\|_{L_\infty (0, t;L_2( \Omega))}^2   \leq   C+C\|u\|_{L_\infty
(0, t;L_2( \Omega))},$$ so $$  \|u\|_{L_\infty (0, t;L_2(
\Omega))}   \leq   C,$$ and (5.4) follows from (5.15). $ \Box $

Estimate (5.4) means, in particular, that in a certain sense the concentration $u(t)$ tends to a constant as $t\to\infty$.

\end{document}